\documentclass[12pt]{amsart}
\usepackage{amscd,amssymb,latexsym}
\theoremstyle{plain}
\newtheorem{Rmk}[equation]{Remark}
\newtheorem{Thm}[equation]{Theorem}

\newtheorem{Cor}[equation]{Corollary}
\newtheorem{Prop}[equation]{Proposition}
\newtheorem{Lem}[equation]{Lemma}
\newtheorem{Def}[equation]{Definition}
\newtheorem{Fact}[equation]{Fact}
\numberwithin{equation}{section}

\def\bt{\begin{Thm}}
\def\et{\end{Thm}}
\def\br{\begin{Rmk}}
\def\er{\end{Rmk}}
\def\bc{\begin{Cor}}
\def\ec{\end{Cor}}
\def\bp{\begin{Prop}}
\def\ep{\end{Prop}}
\def\bl{\begin{Lem}}
\def\el{\end{Lem}}
\def\bd{\begin{Def}}
\def\ed{\end{Def}}
\def\bq{\begin{quotation}}
\def\eq{\end{quotation}}
\def\beq{\begin{equation}}
\def\eeq{\end{equation}}
\def\bea{\begin{eqnarray}}
\def\eea{\end{eqnarray}}
\def\bfa{\begin{Fact}}
\def\efa{\end{Fact}}
\def\bex{\begin{eqnarray*}}
\def\eex{\end{eqnarray*}}

\newcommand{\mc}[1]{{}}

\newcommand{\m}[1]{\mathbb{ #1}}
\renewcommand{\c}[1]{\mathcal{ #1}}

\def\al{\alpha}               
         
\def\eta{\eta}         \def\th{\theta}       
             
\def\si{\sigma}                
\def\ph{\varphi}               
\def\om{\omega}              
             \def\Ph{\Phi}

\def\ra{\rightarrow}
\def\hb{\hfill$\Box$}

\def\noi{\noindent}

\def\vs{\vspace{1em}}

\def\proof{\noi {\bf Proof }}

\def\GL{{\rm GL}}

\def\and{{\rm and}}

\renewcommand{\O}{\rm{O}}

\renewcommand{\br}{\m R}
\begin{document}

\title[Polar decomposition]{Polar decomposition for p-adic symmetric spaces}

\author{Yves Benoist and Hee Oh}

\address{DMA-ENS 45 rue d'Ulm Paris 75005}

\thanks{the second author is partially supported by NSF grant 0629322}

\email{benoist@dma.ens.fr}

\address{Box 1917, 151 Thayer street, Providence, RI}

\email{heeoh@math.brown.edu}
\begin{abstract}{Let $G$ be the group of
$k$-points of a connected reductive $k$-group
and $H$ a symmetric subgroup associated to an involution $\sigma$ of $G$.
We prove a polar decomposition $G=KAH$ for the symmetric
space $G/H$ over any local field $k$ of characteristic not $2$.
Here $K$ is a compact subset of $G$ and $A$ is a finite union of
groups $A_i$ which are the $k$-points of
maximal $(k,\si)$-split tori, one for each $H$-conjugacy class.

This decomposition is analogous to the well-known
polar decomposition $G=KAH$ for a real symmetric space $G/H$.}
\end{abstract}
\maketitle

%\centerline{Preliminary draft}
%\sidenote{changed}

%1
\section{Introduction}
\label{secintro}

%\subsection{General setting}
%\label{secgeneral}

We begin with a short motivation of this article.
 It is well known that {\it every non-degenerate quadratic form
on $\m R^n$ can be put into diagonal form by an orthogonal
base change}.

This statement is an instance of the so-called polar decomposition
$G=KAH$ for a real symmetric space $G/H$
(see section \ref{secreal}), with
$G={\GL}(n,\m R)$, $H={\O}(p,n\! -\! p)$, $1\le p\le n$,
$K={\O}(n)$ and $A$ the subgroup of diagonal matrices of $G$.

Let $k$ be a non-archimedean local  field.
% and ${\cal O}$ its ring of integer.
It is also well known that
 {\it if the residual characteristic of $k$ is not $2$,
every non-degenerate quadratic form
on $k^n$ can be put into diagonal form by a base change
preserving the sup norm} (see section \ref{secquadratic}).
The aim of this paper is to prove a polar decomposition
for any symmetric space $G/H$ over a local field of characteristic not $2$

Here is the precise statement of our polar decomposition.
Let $k$ be a non-archimedean local field of characteristic not $2$,
${\bf G}$ a connected reductive $k$-group, $\si$
a $k$-involution of ${\bf G}$ and ${\bf H}$ an open
$k$-subgroup of the group ${\bf G}^\si$ of $\si$-fixed points.
A  $k$-torus ${\bf S}$ of ${\bf G}$
is said to be
$(k,\si)$-split if it is $k$-split
and $\si(g)=g^{-1}$ for all $g\in {\bf S}$.

We use the notation $G$, $H$, $S$,... to denote the groups of $k$-points of
${\bf G}$, ${\bf H}$, ${\bf S}$,...

By a theorem of Helminck and
Wang [HW], there are only finitely many $H$-conjugacy classes
of maximal $(k, \sigma)$-split tori of ${\bf G}$.
Choose a set $\{{\bf A}_i:1\le i\le n\}$ of representatives of
$H$-conjugacy class of maximal $(k,\si)$-split tori
of ${\bf G}$ and set $A=\cup_{i=1}^n A_i$.

\bt
\label{thkah}
%{\bf (polar decomposition of symmetric spaces)}\\
There exists a compact subset $K$ of $G$ such that we have $$
G=KAH\; .
$$
\et

We remark that all the subgroups $A_i$, $1\le i\le n$ are needed,
up to $H$-conjugacy, in the above decomposition
(Lemma \ref{lemuseful}) and that we can not
in general choose $K$ to be a compact subgroup of $G$.

The proof relies simultaneously on algebraic properties
of the symmetric spaces $G/H$ mostly due to
Helminck and Wang,
as well as on geometric properties of the
Bruhat-Tits building $X_G$ of $G$.
In fact, we will deduce Theorem \ref{thkah}
from the following geometric statement:
 (see (\ref{defsiflat}) for the definition of $\sigma$-apartment):
\bt There exists a constant $C>0$ such that
any point $x\in X_G$ has distance at most $C$ from
some $\sigma$-apartment in $X_G$.
\et

This paper is organized as follows:
 in section \ref{secbruhattits} we recall
the construction of the Bruhat-Tits building of $G$,
as well as the polar decomposition for the group $G$ itself
which is a consequence of the Bruhat-Tits theory.
After recalling in section \ref{sechelminckwang} a few known
facts on
symmetric spaces over local fields,
we prove in section \ref{seckah} a geometric version of the
polar decomposition  (Theorem \ref{thkah2}) from which we deduce
Theorem \ref{thkah}.
We finish in section \ref{secexample},
by giving a few examples and comments. It may be a good idea
for a reader to start with
this section which is independent of the rest of the article.
\vs

We plan to apply our polar decomposition given in
Theorem \ref{thkah} to some counting
and equidistribution results,
with rate of convergence,
for $S$-integral points on symmetric spaces
in a forthcoming paper [BO].
That is why we  do not work here in the more general setting of
a valued field.
\vs

\noi{\footnotesize
After a lecture by the first author in CIRM
on this paper, P. Delorme
gave him the preprint \cite{DS} with V. S\'echerre,
written simultaneously and independantly,
in which they prove
similar results when the residual characteristic is not $2$.
Their method uses the Bruhat-Tits buildings
in a more algebraic way.}

\section{Bruhat-Tits buildings and polar decomposition}
\label{secbruhattits}

Let $k$ be a non-archimedean local field with a valuation $\om$
and $G$ be the group of $k$-points of a reductive $k$-group ${\bf G}$.
In this section we recall the definition
of the Bruhat-Tits building $X_G$ of $G$, which is
 a metric space on which
$G$ acts properly by isometries
with a compact fundamental domain. This metric space is very
similar to the Riemannian symmetric space associated to
a real reductive group. It is a very important tool
since it gives a geometric insight in understanding the algebraic
properties of $G$.
It is precisely this insight that we will use to prove our main theorem
in this paper.
We also recall several properties of $X_G$ that we will need later.

The construction of $X_G$ relies on a cautious construction
of certain open compact subgroups of $G$.
All the theory on building presented below
 is due to Bruhat and Tits (\cite{BT1}, \cite{BT2}) but the readers may find
 references such as  \cite{Ti}, \cite{Ro2} or \cite{La}
shorter and more helpful.

%21
\subsection{Euclidean Buildings}
\label{secbuilding}
We first give the axiomatic definition of a building.
All our buildings will be Euclidean and discrete.

Let $E$ be a real affine space endowed with
an Euclidean distance $d_E$ and
let $W$ be a discrete subgroup, generated by (hyperplane) reflections,
of the group of affine isometries of $E$.
Let $V$ be the corresponding vector space
and $^vW\subset {\rm GL}(V)$
the finite group given by the linear part of $W$.
Here are a few vocabularies.
A {\it wall} is a hyperplane of $E$ pointwise fixed
by a reflection of $W$.
A {\it root} is a halfspace of $E$ bounded by a wall.
A {\it facet} is an equivalence class in $E$ with respect to the relation:
$x\sim y$ if and only if $x$ and $y$ live in the same roots.
The {\it type} of a facet is its $W$-orbit.
A {\it chamber} is an open facet.
A point $x$ in $E$ is {\it special}
if every wall is parallel to a wall containing $x$. Such a point
always exists.
\vs

Buildings are always associated
to such a pair $(E,W)$. Roughly, they are obtained
by gluing together copies of $E$ along convex union of facets
in such a way that they satisfy very strong geometric properties.
More precisely:
\vs

A {\it building}, modeled on $(E,W)$, is a metric space $(X,d)$
endowed with a {\it system of apartments},
i.e., a collection of subsets $\c A$ of $X$ called {\it apartments},
and a partition of $X$ into subsets called {\it facets}
such that :\\
- for each apartment $\c A$ there is an isometry from $E$ to $\c A$
sending facets to facets,\\
- any two points of $X$ is contained in at least one apartment,\\
- for any apartments $\c A $ and $\c A'$, the intersection $\c A\cap \c A'$
is a closed convex subset of $\c A$ which is a union of facets and
there is a facet preserving isometry from $\c A$ to $\c A'$ which is
the identity on $\c A\cap \c A'$.
\vs

As a metric space, every building $X$ is CAT$(0)$-space (see \cite{KL}).
A {\it translation} of $X$ is an isometry of $X$ which
induces a translation in each apartment $\c A\simeq E$ of $X$.
The building $X$ is {\it essential} if the only translation
of $X$ is the identity.
Denote by $V_0$ the vector space of translations of
$X$. The quotient
$X^{\rm ess}:=V_0\backslash X$ is naturally an essential building
which is called the {\it essential quotient} of $X$
and $X$ is isometric to the Euclidean product $X^{\rm ess}\times V_0$.

%22
\subsection{Reductive k-groups}
\label{secreductive}
We now recall the construction of some important
compact open subgroups of $G$
(which are called, at least when ${\bf G}$
is semisimple, connected and simply-connected,
{\it parahoric subgroups}).

For any $k$-split torus ${\bf S}$ of ${\bf G}$, let
$X^*(S)$ be the group of $k$-characters of ${\bf S}$,
$X_*(S)$  the group of $k$-co-characters of ${\bf S}$,
$V(S)$ the real vector space
$V(S)=X_*(S)\otimes_{\m Z}\m R$.
Let ${\bf T}$ be a maximal $k$-split torus  of ${\bf G}$
so that $\dim ({\bf T})=k$-${\rm rank}({\bf G})$.
Let ${\bf N}$ and ${\bf Z}$ be respectively
the normalizer and centralizer of ${\bf T}$ in ${\bf G}$,
and $\Phi$ be the root system of ${\bf G}$ relative to ${\bf T}$.
For any root $\al\in \Ph$, let ${\bf U}_\al$ be
the corresponding root subgroup.
The group $G$ is generated by $N$ and these root groups $U_\al$ \cite{Bo}.
For any $u$ in $U_\al$, $u\neq 1$, there exists
a unique element
$$m(u)\in N\cap U_{-\al}uU_{-\al}\; .
$$
The image of $m(u)$ in the group $N/Z$
 is the reflection $s_\al$ associated with
the root $\al$.
We denote by $^vW$ the group generated by
 these reflections $s_\alpha$, $\alpha\in \Phi$,
and call it the Weyl group (associated to ${\bf T}$).
When ${\bf G}$ is connected one has $^vW=N/Z$.

The normalizer $N$ acts linearly on the real vector space $V(T)$.
Choose an $N$-invariant euclidean structure on $V( T)$.
There is a unique morphism $\nu$ of the
group  $Z$ into the vector space
$ V(T)$
given by the formula
$$\chi(\nu(z))=-\om(\chi(z))\; ,
$$
for all $z\in Z$ and $\chi\in X^*(Z)$ where
on the left-hand
side $\chi$ is seen as a linear form on $V(T)$.
Hence $Z$ acts by translation via $\nu$ on
the affine space  $E$ underlying $V(T)$.
Since $N/Z$ is a finite group, one has $H^2(N/Z,V(T))=0$
and hence,
this action extends to a morphism
$\nu: N\ra {\rm Aff}(E)$. Such an extension is unique up to
a translation.
The hyperplane $H(u)$ fixed by $\nu(m(u))$
is defined by the equation $\al (x) +\ph_\al(u)=0$ for some
$\ph_\al(u)\in \m R$.
For $x\in E$, we introduce the set
$$
U_{\al,x}:=\{ u\in U_\al\; \mid \; u=1
\;\; {\rm or}\;\;
\al(x)+\ph_\al(u)\geq 0\}\; .
$$
By \cite{BT2}, $U_{\al, x}$ is in fact a subgroup.
We also set
$
N_x:=\{ n\in N \; \mid \; \nu(n)(x)=x\}\;
$, which is clearly a subgroup.

Let $K_x$ denote the subgroup of $G$ generated by $N_x$ and $U_{\al,x}$,
 $\al\in \Ph$. The subgroup
$K_x$ is a compact open subgroup of $G$.

%23
\subsection{Bruhat-Tits Buildings}
\label{secbtb}
The Bruhat-Tits building $X_G$ of $G$
is now defined to be the quotient of $G\times E$
by the equivalence relation :
$$
(g,x)\sim (h,y) \Longleftrightarrow
\exists n\in N\text{ such that }  y=\nu(n)x \;\;{\rm and}\;\;  g^{-1}hn\in K_x \; .
$$
This building is modeled on $(E,W)$ where
$W$ is the {\it affine Weyl group} i.e.
the group generated by the affine reflections
$\nu(m(u))$ with $u\in U_\al$, $u\neq 1$  and $\al\in \Ph$.

The apartments of $X_G$ are parametrized by
 maximal $k$-split tori of $G$ :
the apartment $\c A_{T'}$ of $X_G$ corresponding to the
torus ${\bf T}'=g{\bf T}g^{-1}$ for some $g\in G$ is the image
 of $g\times E$ in $X_G$.

The construction of $X_G$ does not depend on the choice of
a maximal $k$-split torus ${\bf T}$.

The $G$-action by left translations
on the first factor of the product $G\times E$
induces an isometric action of $G$ on $X_G$, which is
 proper and co-compact.
The stabilizers in $G$ of points in $X_G$ are
the conjugates of the subgroups $K_x$, $x\in E$.
Moreover, the kernel of this action is
precisely the maximal normal compact subgroup of $G$.

According to (\cite{Ti} 2.1),
the Bruhat-Tits building can be characterized as
in the following proposition:

\bp \label{probuilding}
{\bf \cite{Ti} }
The Bruhat-Tits building $X_G$ of $G$ is
the unique $G$-set
containing a subset $\c A$ normalized by $N$ such that\\
- as an $N$-set, $\c A$ is isomorphic to $E$,\\
- every $G$-orbit in $X$ meets $\c A$,\\
- for each $x\in \c A\simeq E$,
the stabilizer of a point $x$ contains $K_x$.
\ep

Uniqueness means that any $G$-set satisfying the above
three properties is
$G$-isomorphic to $X_G$. The $G$-isomorphism
is unique only up to a translation of $X_G$.
The vector space $V_0$ of translations of $X_G$ is
equal to $V(S_0)$ where
${\bf S_0}$ is the maximal $k$-split torus
of the center of ${\bf G}$.

%24
\subsection{Polar decomposition for $G$}
\label{secpolar}
We now recall
the polar decomposition (also called Cartan decomposition)
for the group $G$.
This decomposition is an algebraic
corollary of the following geometric fact :
$G$ acts strongly transitively on the building $X_G$ i.e.,
 $G$ acts
transitively on the set of all pairs $(C,\c A)$ where $C$ is a chamber
and $\c A$ is an apartment containing $C$.

Choose a positive root system $\Ph^+$ in $\Ph$ and set
$$Z^+:=\{ z\in Z\;\mid \; \om(\al(z))\leq 0\; ,\; \forall\al\in
\Ph^+\} .$$

We also choose a special
point $x$ contained in the apartment $\c A_T$.

\bp
\label{probt}
{\bf \cite{BT1} }
For any non-archimedean local field $k$ and
 any reductive $k$-group ${\bf G}$,
we have $$G=K_xZ^+K_x. $$
\ep

The special case of the above decomposition
when $G={\rm GL}(n,k)$
 has the following interpretation,
when we consider the action of $G$ on the set of
all ultra-metric norms
on $k^n$.

Recall that a norm $N$ on $k^n$ is called {\it ultra-metric} if
$$N(x+y)\leq {\rm max}(N(x),N(y))\quad
\text{ for all $x$, $y\in k^n$}.$$
Noting that $K_x$ is any conjugate of
${\rm GL}(n, \mathcal O)$ for
the valuation ring $\mathcal O$ of $k$ and
that $Z^+$ is a conjugate of the semi-group consisting of
all diagonal matrices in valuation decreasing order,
Proposition \ref{probt} for ${\rm GL}(n, k)$ implies
that {\it  for any two ultra-metric norms $N_1$ and $N_2$ on $k^n$
there exists a basis of $k^n$ with respect to
 which both of $N_1$ and $N_2$ are diagonal
i.e. of the form
$N(x_1,\ldots ,x_n)=
\sup_{1\leq i\leq n}(\al_i|x_i|)$ with
$\al_i>0$.}

%3
\section{Structure of symmetric spaces}
\label{sechelminckwang}
We collect in this section
 a few preliminary facts on symmetric spaces
due to Helminck and Wang (see \cite{HW} and also \cite{HH}).
\vs

Let $k$ be a field with ${\rm char}(k)\neq 2$,
${\bf G}$ a connected reductive $k$-group, $\si$
a $k$-involution of ${\bf G}$ and ${\bf H}$ an open
$k$-subgroup of ${\bf G}^\si$.

%31
\subsection{Existence of $(k,\si)$-split tori}
\label{seckrank}
Recall that a $k$-torus ${\bf S}$ of ${\bf G}$
is said to be
$(k,\si)$-split if it is $k$-split
and $\si(g)=g^{-1}$ for all $g$ in ${\bf S}$.

Here is a criterion for the existence
of a $(k,\si)$-split torus in ${\bf G}$ (Proposition 4.3 of \cite{HW}).

\bl
\label{lemrank0}
{\bf \cite{HW} }
The following are equivalent:
\begin{enumerate}
\item[(i)] Every $(k,\si)$-split torus of ${\bf G}$ is trivial.
\item[(ii)] Every $k$-split torus of ${\bf G}$ is pointwise fixed by $\si$.
\item[(iii)] Every minimal normal isotropic $k$-subgroup of ${\bf G}$
is pointwise fixed by $\si$.
\end{enumerate}
\el

%32
\subsection{Finiteness of $(k,\si)$-split tori}
\label{secksisplittori}
We will not need the whole description of
the $H$-conjugacy classes of
maximal $(k,\si)$-split tori of ${\bf G}$.

We will only need the following proposition, which can be found in
(4.5), (10.3), (6.10) and (6.16) of \cite{HW}.

\bp
\label{prohw1}
{\bf \cite{HW} }
\begin{enumerate}
\item [(a)] Any maximal $k$-split torus  of ${\bf G}$
containing a maximal $(k,\si)$-split torus
of ${\bf G}$ is $\si$-invariant.
\item [(b)] Any two maximal $(k,\si)$-split tori of ${\bf G}$ are
conjugate to each other by an element of $G$.
\item [(c)] If $k$ is a local field, the number of $H$-conjugacy classes
of maximal $(k,\si)$-split tori of ${\bf G}$ is finite.
\end{enumerate}
\ep

Hence the maximal $(k,\si)$-split tori of ${\bf G}$
have the same dimension, which is defined to be the
$k$-rank of $G/H$.

%33
\subsection{Finiteness of $\si$-invariant $k$-split tori}
\label{secksplit}
Once again, we will not need the precise description of
the $H$-conjugacy classes of
$\si$-invariant maximal $k$-split tori of ${\bf G}$,
but only need the following proposition.

The claim (a) below is a straightforward consequence of
Proposition 2.3  of \cite{HW}
applied to the centralizer of ${\bf T}$.
For other claims,
we refer to
Proposition 2.3, Lemma 2.4 and Corollary 6.16 of \cite{HW}.
We recall that every maximal $k$-split torus of
${\bf G}$ is contained in a minimal parabolic
$k$-subgroup ${\bf P}$ of ${\bf G}$.

\bp
\label{prohw2}
{\bf \cite{HW} }
\begin{enumerate}
\item [(a)] Any  maximal $\si$-invariant $k$-split torus
${\bf T}$ of ${\bf G}$
is a maximal $k$-split torus of ${\bf G}$.
\item[(b)] Any minimal parabolic
$k$-subgroup ${\bf P}$ of ${\bf G}$ contains a
maximal $\si$-invariant $k$-split torus of ${\bf G}$,
unique up to conjugacy by $H\cap R_u(P)$.
\item[(c)] If $k$ is a local field, the number of $H$-conjugacy classes
of minimal  parabolic $k$-subgroups is finite.
Hence the number of $H$-conjugacy classes of
maximal $\si$-invariant $k$-split tori of ${\bf G}$ is finite as well.
\end{enumerate}
\ep

When $k$ is local with ${\rm char}(k)=0$,
the claims (c) in both Propositions \ref{prohw1} and \ref{prohw2}
are also consequences of the following theorem of Borel and Serre:
for every $k$-group ${\bf L}$ and its $k$-subgroup ${\bf M}$,
the number of $L$-orbits in the $k$-points
of  ${\bf L}/{\bf M}$ is finite
(th. 5 of III.4.4 in \cite{Se}).

%4
\section{$KAH$ decomposition}
 \label{seckah}
To prove in Theorem \ref{thkah},
we first interpret it in geometric terms which
is possible, using the
Bruhat-Tits building.

%41
\subsection{Flats, parallelism and Levi subgroups}
\label{secflat}
We need a few more results on the geometry
of Euclidean buildings.
\vs

Let $X$ be an Euclidean building.
A {\it flat} in $X$ is a subset isometric to some Euclidean vector
space. A {\it geodesic} is a one-dimensional flat.
The system of apartments of $X$ is {\it complete}
if it is maximal or, equivalently,
if every maximal flat is an apartment.
In fact, every system of apartments can be extended to a complete
system of apartments.
Moreover, the system of apartments is complete
if and only if every geodesic of $X$ is contained in some
apartment (\cite{Pa} Prop. 2.18).
In this case every flat in $X$ is contained in some apartment.

By (\cite{BT1} 2.8.4), the  system of apartments of
the Bruhat-Tits building $X_G$
is complete.

We now suppose that the system of apartments of $X$ is complete.
Two flats $\c F$ and $\c F'$  of $X$ are called {\it parallel} if their
Hausdorff distance is finite i.e. if there exists
a constant $C>0$ such that
$d(y,\c F')\leq C$ and $d(y',\c F)\leq C$ for every $y\in \c F$,
$y'\in \c F'$.
Parallelism is an equivalence relation among flats.
According to (\cite{KL} 2.3.3 and 4.8.1)~, we have \\
- Two flats $\c F$ and $\c F'$ are parallel if and only if
they are both contained in the same apartment $\c A$
and are parallel in $\c A\simeq E$.\\
- For a given flat $\c F$, the union of
flats $\c F'$ parallel to $\c F$, which is then
the union of apartments containing $\c F$, is
a sub-building of $X$, denoted by $X_{\c F}$.

The {\it singular hull} of a flat $\c F$ is the largest flat $\c F^s$
contained in every apartment containing $\c F$.
A flat is {\it singular} if it is equal to its
singular hull.
By construction, we have $X_\c F=X_{\c F^s}$.

For the Bruhat-Tits building $X=X_G$
of a reductive $k$-group ${\bf G}$, one can describe the
sub-building $X_\c F$ as the Bruhat-Tits building $X_L$
associated to some Levi subgroup ${\bf L}$ of ${\bf G}$
(cf. \cite{Ro1} 3.10), which we explain below.
A {\it Levi $k$-subgroup} ${\bf L}$ of ${\bf G}$
is the centralizer ${\bf Z}({\bf S})$ of some $k$-split torus ${\bf S}$
of ${\bf G}$.
A $k$-split torus ${\bf S}$ in ${\bf G}$ is called
{\it singular} if it is obtained as intersection
of maximal $k$-split tori
or, equivalently,
if it is equal to the maximal $k$-split torus
of the center of some Levi subgroup of ${\bf G}$.
This gives a bijection between the Levi $k$-subgroups
of ${\bf G}$ and the singular $k$-split tori of ${\bf G}$.
The {\it singular hull} of a
$k$-split torus ${\bf S}$ is the
smallest singular $k$-split torus containing ${\bf S}$.

Let ${\bf L}$ be a Levi $k$-subgroup of ${\bf G}$.
Note that any maximal $k$-split torus of ${\bf L}$
is also a maximal $k$-split torus of ${\bf G}$.
According to (\cite{BT2} 4.2.17), the union of the apartments
$\c A_{T}$ of $X_G$ associated to a maximal $k$-split
tori ${\bf T}$ of ${\bf L}$ is a sub-building
isometric to the Bruhat-Tits building
$X_L$ of $L$. This sub-building is also denoted by $X_L$
by abuse of notation.
One can also  check this claim using the
characterization (\ref{probuilding}) of the Bruhat-Tits buildings.

We can summarize the above discussion in the following proposition.

\bp
\label{prosingular}
{\bf (\cite{BT2}, \cite{KL}, \cite{Ro1})}
Let $k$ be a non-archimedean local field and ${\bf G}$ a reductive $k$-group.
\begin{enumerate}
\item [(a)]  For each singular flat $\c F$ in $X_G$,
there exists a unique singular $k$-split torus ${\bf S}$ of ${\bf G}$ such
that $X_\c F= X_{Z(S)}$.
This induces a bijection between
the set
of the parallelism classes of singular flats of $X_G$ and the set
of all singular $k$-split tori of ${\bf G}$.
\item [(b)] The  parallel singular flats associated
to a singular $k$-split torus ${\bf S}$ are exactly the
$V(S)$-orbits in $X_{Z(S)}$.
\end{enumerate}
\ep
\noi {\bf Remarks } - Every apartment
 of $X_{Z(S)}$ is associated to a maximal $k$-split
torus of the Levi subgroup ${\bf Z}({\bf S})$.
Not all of them contain $\c F$.

- The group $Z(S)$ almost never acts transitively on
the set $X_{Z(S)}/V(S)$ of flats parallel to $\c F$,
but it acts co-compactly on it, since it already acts co-compactly
on the building $X_{Z(S)}$.

%42
\subsection{$\si$-flats}
\label{secsiflat}
For the rest of section \ref{seckah},
we suppose as in our initial setting
 that
 $k$ is a non-archimedean local field with ${\rm char}(k)\neq 2$,
${\bf G}$ is a connected reductive $k$-group, $\si$
is a $k$-involution of ${\bf G}$ and that ${\bf H}$ is an open
$k$-subgroup of ${\bf G}^\si$.

We first remark that the involution $\si$ induces an involution,
also denoted by $\si$, of the Bruhat-Tits building $X_G$.
In fact, let
${\bf G'}=\langle \si\rangle \ltimes {\bf G}$
be the semi-direct product of the group $\langle \si\rangle$ of order $2$
with ${\bf G}$. This
is a reductive $k$-group whose Bruhat-Tits  building $X_{G'}$
is the same as $X_G$. That is why $G'$ acts on $X_G$.

The following definition allows us to restate
Theorem \ref{thkah} in a geometric way.

\bd
\label{defsiflat}
A $\si$-flat in $X_G$ is a $\si$-invariant flat on which
$\si$ has only one fixed point.
 A $\si$-flat in $X_G$ of maximal dimension is called
a $\si$-apartment.
\ed

We will see below that the maximal dimension
of a $\sigma$-flat in $X_G$ is equal to the $k$-rank of $G/H$.
Note that a maximal $\si$-flat is not always a $\si$-apartment:
for instance, if $k${\rm -rank}$(G)=k${\rm -rank}$(H)$,
there exists a chamber in $X_G$
of which every point is $\sigma$-fixed. Certainly
each point in this chamber is a maximal $\sigma$-flat
but not of maximal dimension unless
the $k$-rank of $G/H$ is zero.

This shows that
the union of all $\si$-apartments of $X_G$ is not
equal to $X_G$ in general.
However, we will see in Theorem \ref{thkah2} that
this is not far from being true.
\vs

We will need the following two lemmas.
The first one
 is analogous to Lemma \ref{lemrank0}

\bl
\label{lemrank0bis}
\begin{enumerate}
\item [(a)] The following
are equivalent:
\begin{enumerate} \item[(i)] The $k${\rm -rank} of $G/H$ is
equal to $0$.
\item[(ii)] Every $\si$-apartment of $X_G$ is a
point.
\item[(iii)] The involution $\si$ acts trivially on $X_G$.
\end{enumerate}
\item[(b)] Any point $x\in X_G$ with $\si(x)\neq x$ lies in a
one-dimensional $\si$-flat.
\end{enumerate}
\el

\proof
$(i)\Rightarrow (ii)$
Let $\c F$ be a $\si$-flat of dimension $1$.
Note that the singular hull $\c F^s$ of $\c F$ is a $\si$-invariant
flat of $X_G$. By Proposition \ref{prosingular},
there exists a singular $k$-split torus ${\bf S}$ such that
$X_{\c F^s}=X_{Z(S)}$. This torus ${\bf S}$
and its centralizer ${\bf L}$ is $\si$-invariant.
The building $X_L$ is isometric to a product
$X_L^{\rm ess}\times V(S)$.
The group $S$ acts on this product by translations on the Euclidean
factor $V(S)$.
Since the action of $\si$ on $V(S)\simeq \c F^s$ is non-trivial,
the action of $\si$ on $S$ is also non-trivial. It follows
that the $k${\rm -rank} of $G/H$ is non-zero.

$(ii)\Rightarrow (iii)$
For $x\in X_G$, suppose $\si(x)\neq x$.
We want to prove that there exists a one dimensional $\si$-flat
$\c F$ containing $x$. This will also prove the claim (b).
Since there exists an apartment containing both $x$ and $\si(x)$,
there exists a geodesic $\c F_0$ containing both $x$ and $\si(x)$.
The points $x$ and $\sigma(x)$
 cut this geodesic into three pieces :
two rays $\c F_-$ and $\c F_+$ and a segment $[x,\si(x)]$.
Let
$$
\c F:= \c F_-\cup [x,\si(x)] \cup \si(\c F_-)\; .
$$
This path $\c F$ is indeed
a geodesic since it is locally a geodesic
and $X_G$ is a CAT$(0)$-space. By construction, $\c F$ is
also $\si$-invariant.
Since $\si(x)\neq x$, $\c F$ is a one-dimensional $\si$-flat,
as desired.

$(iii)\Rightarrow (i)$
Let ${\bf S}$ be a maximal $(k,\si)$-split
torus of ${\bf G}$ and ${\bf T}$
a maximal $k$-split torus of ${\bf G}$
containing ${\bf S}$.
By Proposition \ref{prohw1} (a), the torus
${\bf T}$ is $\si$-invariant. Hence the
corresponding apartment $\c A_T$ is also $\si$-invariant.
If the $k$-rank of $G/H$ is non-zero, the action of $\si$ on
${\bf T}$ is non-trivial. Hence the action of $\si$ on
$\c A_T$ is also non-trivial. This proves the claim.
\hb

%43
\subsection{$\si$-apartments}
\label{secsiapartment}
Note that a maximal $(k,\si )$-split torus of ${\bf G}$ may not be
a singular $k$-split torus.
Nevertheless, we have the following assertion
analogous to Proposition \ref{prosingular}.

\bl
\label{lemsiapartment}
\begin{enumerate}
\item[(a)] For any $\sigma$-apartment $\c F$ of $X_G$,
there exists a unique
maximal $(k,\si )$-split tori ${\bf R}$ of ${\bf G}$ such that
  $X_\c F=X_{Z(R)}$.
This induces a bijection between the
set of
all parallelism classes of $\si$-apartments of  $X_G$
and
the set
of all maximal $(k,\si )$-split tori of ${\bf G}$

\item[(b)] For a maximal $(k,\si )$-split torus ${\bf R}$ of ${\bf
G}$, the parallel $\si$-apartments $\c F$ associated to ${\bf R}$
are exactly the
$V(R)$-orbits in the building $X_{Z(R)}$ and
the subgroup $H\cap Z(R)$ acts co-compactly on the quotient building
$X_{Z(R)}/V(R)$.
\item[(c)] Every apartment of  $X_G$ containing a $\si$-apartment
is $\si$-invariant.
\end{enumerate}
\el

As an immediate corollary, we obtain:
\bc
\label{corsiapartment}
The $k$-rank of $G/H$ is equal to
the dimension of the $\si$-apartments of $X_G$.
\ec

\proof {\bf of Lemma \ref{lemsiapartment} }
We will prove (a), (b) and (c) simultaneously.

Let $\c F$ be a $\si$-apartment of $X_G$.
Let $\c F^s$ be the singular hull of $\c F$ and ${\bf S}$
the unique singular $k$-split torus of ${\bf G}$ such that
$X_{Z(S)}=X_{\c F^s}$ (see Proposition \ref{prosingular}).
Note that the torus ${\bf S}$ and its centralizer
${\bf L}:={\bf Z}({\bf S})$ are also $\si$-invariant.
Let ${\bf R}$ be the maximal $(k,\si)$-split torus of ${\bf S}$.
Again by Proposition \ref{prosingular}, the vector space
$V(S)$ acts simply transitively on $\c F^s$.
Since this action is $\si$-equivariant, and $\c F$ is a
$\si$-flat of maximal dimension
in $\c F^s$,  $\c F$ is a $V(R)$-orbit.

Since $\c F^s$ is the singular hull of $\c F$,  the torus
${\bf S}$ is also the singular hull of ${\bf R}$. It follows that
${\bf L}$ is equal to the centralizer of ${\bf R}$ and
hence $X_{Z(R)}=X_{\c F}$.

By the maximality of the $\si$-flat $\c F$, the quotient building
$X_L/V(R)$ does not contain any one-dimensional $\si$-flat. Hence,
by Lemma \ref{lemrank0bis} (i),
${\bf L}/{\bf R}$ does not contain any non-trivial $(k,\si)$-split
torus and hence
 ${\bf R}$ is a maximal $(k,\si)$-split torus in ${\bf L}$.
Since ${\bf L}$ is the centralizer of ${\bf R}$,
it follows that
${\bf R}$ is a maximal $(k,\si)$-split torus of ${\bf G}$.
To show the uniqueness, if ${\bf R}'$ is any maximal
$(k,\si)$-split torus of ${\bf G}$ such that $X_{Z(R')}=X_\c F$,
then ${\bf R}'$ is contained in ${\bf S}$ and hence in ${\bf R}$ as well.
By the maximality of ${\bf R}'$, ${\bf R'}={\bf R}$.
This shows the first part of the claim (a).

By Lemma \ref{lemrank0} (iii), the group $H\cap L$ acts co-compactly
on $X_L/V(R)$. This, together with the previous discussion, proves (b).
By Lemma \ref{lemrank0bis} (iii), the action of $\si$ on
the building $X_L/V(R)$ is trivial. Therefore every flat containing
$\c F$ is $\si$-invariant. This proves (c).

To finish the proof of (a), letting
 ${\bf R}$ be a maximal $(k,\si)$-split torus  of
${\bf G}$, we want to prove that there exists
a $\si$-apartment $\c F$ of $X_G$, unique
to parallelism, such that $X_{Z(R)}=X_\c F$.

The singular hull ${\bf S }$ of ${\bf R }$ and its centralizer
${\bf L}={\bf Z}({\bf S})$, being
 equal to ${\bf Z}({\bf R})$, are $\si$-invariant.
The associated building $X_L$ is also $\si$-invariant and the action
of the vector space $V(S)$ on $X_L$ commutes with $\si$.
Hence the $V(R)$-orbits in $X_L$ are $\si$-flats
any two of which are parallel to each other.
 Let $\c F$ be one of these $\sigma$-flats.
Since ${\bf R}$ is maximal, the reductive $k$-group
${\bf L}/{\bf R}$ does not contain any non-trivial
$(k,\si)$-split torus.

By Lemma \ref{lemrank0bis} (iii), $\si$ acts trivially on
the quotient building $X_L/V(R)$. Hence $\c F$ is a
$\si$-flat of maximal dimension, i.e., a $\si$-apartment.
Since the torus ${\bf S}$ is the singular hull of ${\bf R}$,
the $V(S)$-orbit $\c F^s$ containing $\c F$ is the singular hull of $\c F$
and we have
$X_{Z(R)}=X_\c F$.

Finally, if $\c F'$  is any $\si$-apartment
such that $X_{Z(R)}=X_{\c F'}$,
then $\c F'$ is contained in a $V(S)$-orbit and  hence in a $V(R)$-orbit
in $X_L$.
Therefore since $\c F'$ is of maximal dimension,
$\c F'$ is equal to a $V(R)$-orbit in $X_L$.
Hence $\c F'$ is parallel to $\c F$.
This completes the proof of (a).

%44
\subsection{Geometric interpretation}
\label{secinterpretation}
Here is a geometric reformulation of Theorem \ref{thkah}.\vs

For $C>0$,
a subset $Y$ of a metric space $X$ is {\it $C$-dense} if
every point of $X$
has distance at most $C$ from $Y$. A subset $Y$
is {\it quasi-dense} in $X$
 if it is $C$-dense for some $C>0$.

\bt
\label{thkah2}
Let $k$ be  a non-archimedean local field with ${\rm char}(k)\neq 2$,
${\bf G}$ a connected reductive $k$-group and
$\si$ a $k$-involution of ${\bf G}$.
The union of all $\si$-apartments
is quasi-dense in the Bruhat-Tits building $X_G$.
\et

We begin by proving the following lemma:
\bl
\label{lemnonfixed}

Either $\sigma$ acts trivially on $X_G$ or
 the set
$\{ x\in X_G\;\mid \; \si(x)\neq x\}$ is quasi-dense in $X_G$.
\el

\noi {\bf Remark }
When the residual characteristic of $k$ is not $2$,
according to Prasad and Yu \cite{PY},
the set $X_G^\si$ of $\si$-fixed points can be identified
with the Bruhat-Tits building $X_H$ of $H$, when
$\sigma$-action is non-trivial on $X_G$.
This gives a proof of Lemma \ref{lemnonfixed} in this case.
We present below a short proof of the above lemma assuming only
that $\text{char}\, (k)\ne 2$.
\vs

\proof {\bf of Lemma \ref{lemnonfixed} }
Suppose that $\sigma$-action on $X_G$ is non-trivial.
We will show that there exists a constant $C>0$ such that,
for every $x\in X_G$, the ball $B(x,C)$ is not pointwise fixed
by $\si$.

Suppose not. Then there exists a sequence $\{x_n\in X_G\}$
such that the balls $B(x_n,n)$ are pointwise fixed by $\si$.
We may assume that for a fixed $x_0\in X_G$,
 $x_n=g_nx_0$ for some $g_n\in G$.
Then the element $\si_n:= g_n\si g_n^{-1}$ is an involutive
automorphism of
$G$ which acts trivially on $B(x_0,n)$.
Hence the sequence $\si_n$ sub-converges to
the identity
in ${\rm Aut}(G)$.
Hence any compact open subgroup of
${\rm Aut}{(G)}$ contains $\si_n$ for all sufficiently large $n$.
On the other hand,
if ${\rm char}(k)=0$, then ${\rm Aut}(G)$ contains a
compact open subgroup which is torsion free.
If the residual characteristic of $k$ is $p\ne 2$,
the group ${\rm Aut}(G)$,
as any group $G'$ of $k$-points
of a $k$-group ${\bf G}'$,
contains a compact open subgroup
which is a pro-$p$-group.
Since a pro-$p$-group has no elements of order $2$,
the order of $\sigma_n$ is $2$, and we obtain a contradiction
in either case.
\hb\vs

\proof {\bf of Theorem \ref{thkah2} }
We want to prove that there exists a constant $C_G>0$ such that
every point $x\in X_G$ has distance at most
$C_G$ from some $\si$-apartment $\c F$ of $X_G$.

We argue by induction on the dimension of ${\bf G}$.
Writing, as in section \ref{secbruhattits},
$X_G\simeq X_G^{\rm ess}\times V(S_0)$ where
${\bf S}_0$ is the maximal $k$-split torus of
the center of ${\bf G}$, we may assume that ${\bf G}$ is
semisimple.

If $\si$ acts trivially on $X_G$ then
every point of $X_G$ is a $\si$-apartment in which case
our claim trivially holds.
Hence we assume that $\si$ does not act trivially on $X_G$.
By Lemma \ref{lemnonfixed} and Lemma \ref{lemrank0bis}
(b), we may then assume
that
$x$ lies in a one-dimensional $\si$-flat $\c F$.
Then $x$ belongs to the sub-building $X_{\c F^s}$ where $\c F^s$ is
the singular hull of $\c F$.
By Proposition \ref{prosingular} and  Lemma
 \ref{lemsiapartment},
the sub-building $X_{\c F^s}$ is equal to the building $X_L$ associated
to some $\si$-invariant Levi-subgroup
${\bf L }\varsubsetneq {\bf G}$ such that
every $\sigma$-apt in  $X_L$ is also an $\sigma$-apt in $X_G$.
We apply our induction hypothesis to the reductive subgroup ${\bf L}$
to conclude.

Note that the distance on the building
$X_G^{\rm ess}$ is unique only up to some
scalar factors, one for each quasi-simple isotropic
normal $k$-subgroup of ${\bf G}$.
The constant $C_G$ will depend on this normalization of those
scalar factors. The same is true for ${\bf L}$.

Each $\si$-invariant Levi subgroup  ${\bf L}$
gives rise to some constant $C_L$.
This constant $C_L$ depends only on the $H$-conjugacy class
of ${\bf L}$.

To obtain a finite upper bound for these $C_L$, we need to know that
{\it there are only finitely many $H$-conjugacy classes of
$\si$-invariant
Levi subgroup of ${\bf G}$}.
Since a $k$-split torus contains only  finitely many
singular $k$-split sub-tori, this assertion is a direct consequence of
the finiteness of the number of $H$-conjugacy classes of
maximal $\si$-invariant $k$-split tori of ${\bf G}$
(Proposition \ref{prohw2}).
\hb

%45
\subsection{Orbits of $(k,\si)$-split tori in $X_G$}
\label{secorbit}
The following corollary is the link between
Theorem \ref{thkah} and Theorem \ref{thkah2}.

As in the introduction,
 we choose maximal $(k, \si)$-split tori
 ${\bf A}_1,\ldots {\bf A}_n$,
which are representatives of $H$-conjugacy classes
of all maximal $(k, \sigma)$-split tori of ${\bf G}$ (cf. Proposition \ref{prohw1}).
Let $A$ be the union of the
groups of $k$-points $A_i$, $1\le i\le n$.

\bc
\label{corkah2}
Let $k$ be  a non-archimedean local field with ${\rm char}(k)\neq 2$,
${\bf G}$ a connected reductive $k$-group,
$\si$ a $k$-involution of ${\bf G}$ and ${\bf H}$ an
open $k$-subgroup of ${\bf G}^\si$.

For any fixed $x_0\in X_G$,
there exists a constant $C>0$ such that
for all $x\in X_G$, there exist $h\in H$ and $1\le i\leq n$
with
 $$d(hx,A_ix_0)\leq C .$$
\ec

In other words, {\it the union of the $H$-orbits
meeting $Ax_0$ is quasi-dense in $X_G$}.
\vs

\proof
  By Theorem \ref{thkah2}, the set of $\sigma$-apartments are $C_0$-dense
in $X_G$ for some $C_0>0$.
Hence $d(x, x_1)<C_0$ for some $x_1$
contained in a $\si$-apartment $\c F$.
Then by Lemma \ref{lemsiapartment},
 $x_1$ is contained the sub-building $X_{Z(S)}$
of the centralizer of a maximal $(k,\si)$-split torus
${\bf S}$ of ${\bf G}$.
By Proposition \ref{prohw1} (c),
${\bf S}$ is conjugate to one of ${\bf A}_i$ by an element of
$H$ and hence for some $h\in H$,
$hx_1\in X_{Z(A_i)}$.

Fix a $\sigma$-apartment $\c F_i$ associated to $A_i$
as in Lemma \ref{lemsiapartment} for each $1\le i\le n$
and let $e_i$ denote the Hausdorff distance between $\c F_i$ and $A_i x_0$.
By Lemma \ref{lemsiapartment},
the group $H\cap Z(A_i)$ acts co-compactly on the
quotient building $X_{Z(A_i)}/V(A_i)$.
In other words,
 there is $d_i>0$ such that
any point in $X_{Z(A_i)}$
 is within distance $d_i$
 from any
 $\sigma$-apartment in $X_{Z(A_i)}$ parallel to
 $\c F_i$, up to translation by
$H\cap Z(A_i)$.

Therefore
$$d(h'hx_1, \c F_i)<d_i\quad\text{ for some $h'\in H\cap Z(A_i)$} .$$

Hence
$$d(h'hx, A_i x_0)\le d(x, x_1)+d(h'hx_1, \c F_i) + d(\c F_i, A_ix_0)\le
C_0 + d_i + e_i.$$
It remains to put $C:=C_0+\max_{1\le i\le n} (d_i+e_i)$.

\hb

\vs

\proof {\bf of Theorem \ref{thkah} }
It is now very easy to conclude.
Fixing a point $x_0\in X_G$,
let $C>0$ be a constant given by Corollary \ref{corkah2}
and $K$ the compact subset of $G$
given by
$$
K:=\{ k\in G\;\mid \; d(x_0,kx_0)\leq C\}\; .
$$
For any $g\in G$, we apply Corollary \ref{corkah2}
to $x=g^{-1}x_0$ to obtain elements $h\in H$,
 $a\in \cup_{i=1}^n A_i$ such that
$$
d(hx,a^{-1}x_0)\leq C\; ,
$$
or, equivalently,
$$
d(x_0,gh^{-1}a^{-1}x_0)\leq C\; .
$$
In other words, we have $g=kah$ for some $k\in K$.
\hb

%5
\section{Examples}
\label{secexample}

In order to put our polar decomposition in perspective,
we now recall some well-known examples.

%51
\subsection{The real case}
\label{secreal}
We first discuss the archimedean case $k=\m R$ and
recall the proof of the following basic fact:

\bfa
\label{facquadraticr}
Every non-degenerate quadratic form $q$ on the Euclidean space
$\m R^n$ can be put into diagonal form by an orthogonal base
change.
\efa

\proof Choose a point $e_1$ on the Euclidean unit
sphere in $\m R^n$ where $q$ achieves the maximum.
 Note that the Euclidean orthogonal of $e_1$
is also orthogonal for $q$. Hence we can use an
induction argument.
\hb

\vs
\noi {\bf Remark } Fact \ref{facquadraticr} is a special case of
the polar decomposition
$G=KAH$ for a real symmetric space $G/H$ with
$G={\rm GL}(n,\m R)$, $K={\rm O}(n)$,
$A$ the subgroup of diagonal matrices and $ H={\rm O}(p,n\! -\! p )$ for any
$0\le p\le n$:
Write $q=q_0\circ g^{-1}$
 with $q_0= x_1^2+\cdots +x_p^2-x_{p+1}^2-\cdots -x_n^2$~
and $g\in G$.
If we write $g=kah\in KAH$, then
the quadratic form $q\circ k=q_0\circ a^{-1}$ is diagonal.
\vs

We now recall more precisely the ``real'' polar decomposition
and its proof.
Let ${\bf G}$ be a connected reductive $\m R$-group, $\si$
an $\m R$-involution of ${\bf G}$ and ${\bf H}$ an open
$\m R$-subgroup of ${\bf G}^\si$.
Recall that an $\br$-involution
$\th$ of ${\bf G}$ is called
{\it Cartan}
if the group ${\bf K}={\bf G}^\th$ of $\th$-fixed points is a maximal
$\m R$-anisotropic subgroup of ${\bf G}$.

The following facts are well-known (see \cite{Sc}):

- There exists a {\it Cartan involution} $\th$
of ${\bf G}$ which commutes  with $\sigma$,
i.e., $\si\th=\th\si$.

- A maximal $(\m R,\si)$-split torus ${\bf A}$ of
the group ${\bf H}':={\bf G}^{\si\th}$
of $\si\th$-fixed points is a maximal $\m R$-split torus of ${\bf H}'$ as
well as a maximal $(\m R,\si)$-split torus of ${\bf G}$.
Moreover, all maximal $(\m R,\si)$-split tori of ${\bf G}$
are $H$-conjugate.

Letting ${\bf K}={\bf G}^\th$ and
${\bf A}$ a maximal $(\m R,\si)$-split torus of
the group ${\bf H}'$,
 the polar decomposition for $G/H$ is given as follows:

\bp
\label{prokahr}
{\bf (see \cite{Sc}) }
We have $G=KAH$ with $K\! =\! {\bf K}_{\m R}$,
$A\! =\! {\bf A}_{\m R}$ and $H\! =\! {\bf H}_{\m R}$.
\ep

We recall here a sketch of proof to emphasize both analogies and
differences with our proof in the non-archimedean case.

\vs
\noindent {\bfseries Sketch of proof }
As in the proof of Fact \ref{facquadraticr},
 the proof uses a minimizing argument.
Let $g\in G$. Fix a base point $x_0=K/K$ in
the Riemannian symmetric space $X_G=G/K$ and and let
$x:=g^{-1}x_0$.
Let $y=h^{-1}x_0$, with $h\in H$, be the nearest point
to $x$ on the totally geodesic sub-manifold $X_H:=H/(H\cap K)$.
The geodesic from $x_0$ to $hx$ is orthogonal
to $X_H$. Hence the point $hx$ is in the
symmetric subspace
$X_{H'}=H'/(H'\cap K)$ ``orthogonal at $x_0$'' to $X_H$.
Therefore we have $x_0=h'hx$ with $h'\in H'$,
or, in other words, $g=kh'h$ with $k\in K$.
Using the Cartan decomposition
$H'=(H\cap K)A(H\cap K)$ for the subgroup $H'$,
we then obtain $g\in KAH$, as desired.
\hb

%52
\subsection{Quadratic forms}
\label{secquadratic}
We now discuss the diagonalization of quadratic forms in
the non-archimedean case.
The following well-known fact is analogous to Fact \ref{facquadraticr}:

\bfa
\label{facquadratick}
Let $k$ be a non-archimedean local  field
of residual characteristic not $2$.

Every non-degenerate quadratic form
on $k^n$ can be put into diagonal form by a base change
preserving the sup norm.
\efa

Recall that the stabilizer of the sup norm
$\|x\|={\rm sup}_{1\leq i\leq n} |x_i|$
is the maximal compact subgroup
$K_0={\rm GL}(n,\c O)\subset {\rm GL}(n,k)$ where ${\c O}$ is the
valuation ring of $k$.

\vs
\proof
Choose a point $e_1$ on the unit ball $\c O^n$ where $|q|$ is maximum.
Note that $e_1$ is primitive, i.e.,
 $ke_1\cap \c O^n=\c Oe_1$~.
 Hence $e_1$ can be completed to a basis, say $(e_1,\ldots ,e_n)$, of $\c O^n$.
Write $q(\sum x_ie_i)=\sum a_{ij}x_ix_j$ where the matrix $b=(a_{ij})$ is
symmetric. One can assume $a_{11}=1$.
By the maximality of $|q(e_1)|=a_{11}=1$,
we have $|q(e_i)|=|a_{i,i}|\leq 1$ for each $i$.
It follows from $|q(e_1+e_i)|\le 1$
that $|2a_{1,i}|\leq 1$ for each $i$. Since the residual
characteristic is not $2$, we have $|a_{1,i}|\leq 1$.
Therefore the $n$-tuple
$(e_1,e'_2,\ldots ,e'_n)$ with $e'_i:= e_i- a_{1,i}e_1$
is also a basis of $\c O^n$ and $e'_2,\ldots e'_n$ are orthogonal
to $e_1$.
We then use an induction argument to conclude.
\hb

We may translate  Fact \ref{facquadratick} in terms of our polar
decomposition. Let $q_0$ be a non-degenerate
diagonal quadratic form on
$k^n$, ${\bf G}:={\bf GL}(n)$, $K_0:={\rm GL}(n,\c O)$, ${\bf A}_0$ the
$k$-group of diagonal matrices, ${\bf H}:={\bf O}(q_0)$, so that
${\bf X}:={\bf G}/{\bf H}$ is the space of non-degenerate
quadratic forms and $X:=G/H$ is the space of quadratic forms
$G$-equivalent to $q_0$. The subset of diagonal quadratic forms is
$X_{\rm diag}:= X\cap {\bf A}_0q_0$. Fact \ref{facquadratick}
gives the equality
$$
X= K_0
X_{\rm diag}\;\;\;\;{\rm or}\;\;\;\; G=K_0({\bf A}_0{\bf H})_k\; .
$$
The set $({\bf A}_0{\bf H})_k$ is a finite union of double
classes $A_0g_iH$, where the set of indices $i$  is given by
$(k^*/k^{*2})^n$.
Letting
 ${\bf A}_i:=g_i^{-1}{\bf A}_0g_i$, $A$ the union of $A_i$'s and
$K$ the  union of $K_0g_i$'s.
we
 obtain the polar decomposition $G=KAH$ as in Theorem \ref{thkah}.
Note here that $K$ is
a compact \underline{subset}.
\vs

\noi {\bf Remark } If the residual characteristic is $2$ but the
characteristic of $k$ is not $2$, the same argument shows that
there exists a compact \underline{subset} $K_0$ of ${\rm GL}(n,k)$
such that
{\it every non-degenerate quadratic form
on $k^n$ can be put into diagonal form by a base change in $K_0$}
i.e., $G=K_0({\bf A}_0{\bf H})_k$. One can easily find an example
where $K_0$ cannot be chosen to be a subgroup.

%53
\subsection{The group case }
\label{secgroupcase}

The following is a well-known corollary of Proposition \ref{probt}.

\bc
\label{corbt}
Let $k$ be a non-archimedean local field,
${\bf G}$ a reductive $k$-group and ${\bf A}$
a maximal $k$-split torus of $\bf G$. Then
there exists a compact subset $K$ of $G$ such that
$G=KAK$.
\ec

\proof We keep the same notation as in Proposition \ref{probt}.
Let $Z_c=K_x\cap Z$ be the maximal compact subgroup of $Z$.
The quotient $Z/Z_c$ is an abelian group of rank $r$ and
the image of $T$ in $Z/Z_c$
is a subgroup of finite index.
\hb

This corollary is also a special case of our
polar decomposition (Theorem \ref{thkah}) for the group
${\bf G}\times {\bf G}$ with the involution
given by $\si(g_1,g_2)=(g_2,g_1)$ for which ${\bf H}$
is the diagonal embedding of ${\bf G}$ into ${\bf G}\times
{\bf G}$.
In this example there exists only one
$H$-conjugacy class of maximal $(k,\si)$-split tori of ${\bf G}\times {\bf G}$.

\vs
\noi{\bf Remark }
Note that even in this corollary one can not always choose
the compact subset $K$ to be a subgroup of $G$
 because the map from $T$
to $Z/Z_c$ is not surjective in general.

%54
\subsection{All $\si$-apartments are useful}
\label{secuseful}

The following lemma with Lemma \ref{lemsiapartment} (a)
shows that
{\it all the representatives ${\bf A}_i$}
of $H$-conjugacy classes
of maximal $(k,\si)$-split tori of ${\bf G}$
{\it are useful
in our polar decomposition} (Theorem \ref{thkah} or \ref{thkah2}).

\bl
\label{lemuseful}
Let $k$ be a non-archimedean local field with ${\rm char}(k)\neq 2$,
${\bf G}$ a connected reductive
$k$-group, $\si$ a $k$-involution of ${\bf G}$ and ${\bf H}$ an
open $k$-subgroup of ${\bf G}^\si$.

Let $\c F, \c F_1, \ldots ,\c F_\ell$ be maximal $\si$-apartments
(see definition \ref{defsiflat})
of $X_G$ such that
$$
\sup _{x\in \c F}\; d(x, \cup_{1\leq i\leq \ell}H.\c F_i)<\infty\; .
$$
Then
$\c F$ is parallel to one of the translates $h\c F_i$ with $h\in H$
and $1\le i\leq \ell$.
\el

\proof
Since there are only finitely many directions of singular flats
in $\c F$, we can
choose a $\si$-invariant geodesic $t\ra c_t$ in $\c F$ with
the same singular hull as $\c F$ and satisfying $\si(c_t)=c_{-t}$
for all $t$.
By hypothesis we can find $C_0>0$, $1\le i\leq \ell$
and sequences $t_n\ra\infty$,
$h_n\in H$ such that
$$
d(c_{t_n},h_n\c F_i)\leq C_0\; .
$$
Applying $\si$, we have
$$
d(c_{-t_n},h_n\c F_i)\leq C_0\; .
$$
Using the fact that
in a CAT$(0)$-space,
the distance function to a convex subset
is a convex function,
we deduce that for all $s\in [-t_n,t_n]$
$$
d(c_{s},h_n\c F_i)\leq C_0\;
$$ and that
$$
d(c_{0},h_n x_i)\leq C_0\;
$$ where $x_i$ denotes the unique point of $\c F_i$
which is $\si$-invariant.
The sequence $h_n$ remains in a compact subset of $G$ and hence
converges
to an element $h\in H$, after passing to
a sub-sequence. We then have that for all $s\in \m R$,
$$
d(c_{s},h\c F_i)\leq C_0\; .
$$
Hence the geodesic $\c G:=\{  c_t\; ,\; t\in \m R\}$
is parallel to a geodesic in $h\c F_i$
(see 2.3.3 in \cite{KL}).
But then $h\c F_i$ is included in the sub-building $X_{\c G}$.
Since $\c G$ and $\c F$ have the same singular hull,
we have  $h\c F_i\subset X_{\c G}= X_{\c F}$. Therefore
Lemma \ref{lemsiapartment} implies that
$h\c F_i$ is parallel to $\c F$.
\hb

%9

\end{document}